\documentclass[journal,12pt,onecolumn,draftclsnofoot,]{IEEEtran}
\IEEEoverridecommandlockouts
\usepackage{cite}
\usepackage{amsmath,amssymb,amsfonts}
\usepackage{algorithmic}
\usepackage{graphicx}
\usepackage{textcomp}
\usepackage{xcolor}
\usepackage{amsmath}
\usepackage{algorithm}
\def\BibTeX{{\rm B\kern-.05em{\sc i\kern-.025em b}\kern-.08em
    T\kern-.1667em\lower.7ex\hbox{E}\kern-.125emX}}

\title{\LARGE \bf
	A Fully Parallel Primal-Dual Algorithm for Centralized and Distributed Optimization
}

\author{S. Sh. Alaviani and A. G. Kelkar
	\thanks{}
	\thanks{The authors are with the Department of Mechanical Engineering, Clemson University, Clemson, SC, 29634 USA emails:
		{\tt\small salavia@clemson.edu, atul@clemson.edu}. 
	}%
	\thanks{
	}%
}

\begin{document}

	\maketitle
	\thispagestyle{empty}
	\pagestyle{empty}

\begin{abstract}

In this paper, a centralized two-block separable optimization is considered for which a fully parallel primal-dual discrete-time algorithm with fixed step size is derived based on monotone operator splitting method. In this algorithm, the primal variables are updated in an \textit{alternating} fashion like Alternating Direction Method of Multipliers (ADMM). However, unlike existing discrete-time algorithms such as Method of Multipliers (MM), ADMM, Bi-Alternating Direction Method of Multipliers (BiADMM), and Primal-Dual Fixed Point (PDFP) algorithms, that all suffer from sequential updates, all primal and dual variables are updated in \textit{parallel} in the sense that to update a variable at each time, updated version of other variable(s) is not required. One of advantages of the proposed algorithm is that its \textit{direct} extension to multi-block optimization is still convergent. Then the method is applied to distributed optimization for which a fully parallel primal-dual distributed algorithm is obtained. Finally, since direct extension of ADMM may diverge for multi-block optimization, a numerical example of a three-block optimization is given for which the direct extension of the proposed algorithm is shown to converge to a solution.

\end{abstract}

\section{Introduction}

In this paper, we consider the following two-block decomposable optimization (and its extension to multi-block optimization) with affine constraint\footnote{For better comparison between ADMM and our proposed algorithm, we use the same notations for (\ref{1}) as in \cite{boydADMM}.}:

\begin{equation}\label{1}
\begin{aligned}
& \underset{x,z}{\text{min}}
& & f(x)+g(z) \\
& \text{subject to}
& & Ax+Bz=c \\
\end{aligned}
\end{equation}
where $x \in \Re^{n}$ and $z \in \Re^{m}$ are decision variables, $f,g$ are convex functions, and $A \in \Re^{p \times n},B \in \Re^{p \times m}, c \in \Re^{p}$ where $\Re$ denotes the set of all real numbers. The problem (\ref{1}) has found applications in many areas of engineering such as signal and image denoising and restoration \cite{image}, compressed sensing \cite{compressedsensing1}-\cite{compressedsensing2}, and channel estimation and coding \cite{coding}. The research challenge is how to approach an optimal solution of (\ref{1}) by translating the original difficult optimization with respect to primal and dual variables into solving easier problems (unconstrained ones). Motivated by the seminal Arrow-Hurwicz-Uzawa primal-dual dynamics \cite{arrowhurwicz}, several researchers have paid attention to develop primal-dual algorithms to solve (\ref{1}) \cite{wangelia2011}-\cite{unifiedframework}. Note that Arrow-Hurwicz-Uzawa dynamics has been used for distributed optimization in \cite{wangelia2011} (see recent Survey \cite{survey2020} for more references). 

\textit{Method of Multipliers}\footnote{Method of Multipliers is also called \textit{Augmented Lagrangian Method} (ALM) \cite{survey2020}.} (MM) \cite{MM1}-\cite{MM3} is a general-purpose iterative solver for constrained optimization (\ref{1}). MM decomposes the original problem into sub-problems that can be solved easily. MM can handle nonsmoothness of objective functions and generic problem constraints and has strong convergence guarantees \cite{qlearnbertsekas}. Although MM has been beneficial, it suffers from sequential updates. Distributed MM has been investigated by some researchers \cite{mmdist1}-\cite{mmdist3}.

Another iterative solver for (\ref{1}) is \textit{Alternating Direction Method of Multipliers} (ADMM) which was proposed in \cite{ADMM1}-\cite{ADMM2} (see \cite{boydADMM} for more references). Similar to MM, ADMM decomposes the original problem into sub-problems that can be solved easily. ADMM also can handle nonsmoothness of objective functions with strong convergence \cite{boydADMM}. Unlike MM, primal variables $x$ and $z$ in ADMM are updated in an \textit{alternating} fashion that decomposes the updates of primal variables into easier sub-problems. Although ADMM has been beneficial, extension of ADMM to multi-block optimization is still a challenging problem since \textit{direct} extension of ADMM to multi-block is \textit{not} necessarily convergent if no further assumption is imposed (see \cite{ADMMdiverges} for details). In this regard, several variants of ADMM have been investigated to overcome this disadvantage such as 3-block ADMM \cite{ADMM3block1}-\cite{ADMM3block2}, Jacobian decomposition of augmented Lagrangian method with proximal terms \cite{ADMMJacobian}, and prediction-correction methods \cite{couple14}. ADMM and its extensions to multi-block optimization suffer from sequential updates. ADMM has been extended to distributed case in \cite{ADMMdistributedorigin1}-\cite{ADMMdistributedorigin2} for estimation on wireless sensor networks, and some authors have paid much attention to develop several distributed versions of ADMM \cite{018}-\cite{adda}, to cite a few. Recently, \textit{Primal-Dual Method of Multipliers}, which can be taken as an extension of ADMM for solving problems over graphs, has been proposed in \cite{priamdualmm1}-\cite{priamdualmm2}.

\textit{Bi-Alternating Direction Method of Multipliers} (BiADMM) \cite{priamdualmm1}-\cite{biadmm} is another iterative solver for (\ref{1}).  BiADMM iteratively minimizes an augmented bi-conjugate function. Unlike ADMM that always involves three updates (two primal variables and one dual variable) per iteration, BiADMM performs either two or three updates per iteration, depending on the functional construction. Similar to MM and ADMM, BiADMM suffers from sequential updates.

Another iterative algorithm to solve (\ref{1}) is \textit{Primal-Dual Fixed Point} (PDFP) algorithm. PDFP algorithm was originally proposed in \cite{fixedpointprimaldualChen} for separable convex optimization, and some researchers have proposed several PDFP algorithms (see \cite{fixedpointprimaldual1}-\cite{fixedpointprimaldual2} and references therein). PDFP algorithm is based on proximal forward-backward splitting and fixed point iterations (see \cite{fixedpointprimaldualChen} for details). PDFP algorithms have been extended to multi-block optimization in \cite{fixedpointprimaldualmultiblock}. Similar to previous algorithms, PDFP algorithms suffer from sequential updates. 


Recently, a unified framework for existing algorithms for centralized separable convex optimization has been proposed in \cite{unifiedframework}, that also requires sequential updates.  All aforementioned algorithms suffer from sequential updates that make them slow in practice. Also applying \textit{forward-backward-forward method} \cite{modifiedforwardbackwardTseng} to (\ref{1}) does \textit{not} yield parallel algorithms (see \cite[Sec. 4]{modifiedforwardbackwardTseng} for details). Recently, a survey paper \cite{survey2020} considers existing distributed primal-dual algorithms.

\textbf{Contribution:} In this paper, we consider the centralized two-block convex optimization (\ref{1}). By Applying \textit{forward-reflected-backward method} recently proposed in \cite{forwardreflectedbackward}, we derive a fully parallel primal-dual discrete-time algorithm with fixed step size, called \textit{Parallel Alternating Direction Primal-Dual} (PADPD) algorithm. In the PADPD algorithm, the two primal variables $x$ and $z$ are updated in an \textit{alternating} fashion like ADMM, but all (primal and dual) variables are updated in \textit{parallel} in the sense that to update a variable at each time, updated version of other variable(s) is not needed. One of main advantages of the proposed algorithm for centralized optimization is that it can be extended \textit{directly} to \textit{any} finite multi-block optimization while preserving its convergence. Then the approach is applied to distributed optimization for which a fully parallel primal-dual distributed algorithm is obtained.

This paper is organized as follows. In Section II, some preliminaries are given. PADPD algorithm and its extension to multi-block optimization are derived in Section III. In Section IV, the method is applied to distributed optimization. Finally, since direct extension of ADMM diverges for a three-block optimization \cite{ADMMdiverges}, a numerical example of the three-block optimization is given for which the direct extension of the proposed algorithm converges. 

\textit{Notations:} $\Re$ denotes the set of all real numbers. $\mathbb{N}$ denotes the set of all natural numbers. $\emptyset$ represents the empty set. $(.)^{T}$ represents the transpose of a matrix or a vector. For any vector $z \in \Re^{n}, \Vert z \Vert_{2}=\sqrt{z^{T}z},$ and for any matrix $Z \in \Re^{n \times n}, \Vert Z \Vert_{2}=\sqrt{\lambda_{max}(Z^{T}Z)}=\sigma_{max}(Z)$ where $\lambda_{max}$ represents the maximum eigenvalue, and $\sigma_{max}$ represents the largest singular value. Sorted in an increasing order with respect to real parts, $\lambda_{2}(Z)$ represents the second eigenvalue of a matrix $Z$. For any matrix $Z \in \Re^{n \times n}$ with $Z=[z_{ij}]$,  $\Vert Z \Vert_{1}= max_{1 \leq j \leq n} \{\sum_{i=1}^{n} \vert z_{ij} \vert \}$ and $\Vert Z \Vert_{\infty}= max_{1 \leq i \leq n} \{\sum_{j=1}^{n} \vert z_{ij} \vert \}$. $col\{x_{1},\hdots,x_{m}\}:=[x_{1}^{T}, \hdots,x_{m}^{T}]^{T}$ denotes the column vector stacked with vectors $x_{1},\hdots,x_{m}$. $\textbf{0}_{n}$ represents the vectors of dimension $n$ whose elements are all zero. 

\section{Preliminaries}

A vector $v \in \Re^{n}$ is said to be a \textit{stochastic vector} when its components $v_{i}, i=1,2,...,n$, are non-negative and their sum is equal to 1; a square $n \times n$ matrix $V$ is said to be a \textit{stochastic matrix} when each row of $V$ is a stochastic vector. A square $n \times n$ matrix $V$ is said to be \textit{doubly stochastic} when both $V$ and $V^{T}$ are stochastic matrices.

Let $\mathcal{H}$ be a real Hilbert space with norm $\Vert . \Vert $ and inner product $\langle .,. \rangle$. A mapping $H:\mathcal{H} \longrightarrow \mathcal{H}$ is said to be $K$\textit{-Lipschitz continuous} if there exists a $K > 0$ such that 
$\Vert Hx-Hy \Vert \leq K \Vert x-y \Vert$ for all $x,y \in \mathcal{H}$.

\textbf{Definition 1:} \textit{Resolvent} of an operator $A$ is defined as $J_{A}:=(I+A)^{-1}$ where $I$ denotes the Identity operator.

\textbf{Remark 1} \cite{minty}-\cite{maximalandproximalRackafler}: If $f: \mathcal{H} \longrightarrow \Re \cup \{+\infty\}$ is a proper closed convex function, the resolvent of its subdifferential, namely $J_{\partial f}$, is equal to the Moreau's proximity operator of $f$ \cite{moreauproximity}, i.e.,
$$prox_{f}(x):=\arg\min_{z \in \mathcal{H}}[f(z)+\frac{1}{2} \Vert z-x \Vert^{2}].$$

\textbf{Definition 2:} A multi-valued mapping $T:\mathcal{H} \rightrightarrows \mathcal{H}$ is called \textit{monotone operator} if $\langle x'-y',x-y \rangle \geq 0$ whenever $x' \in T(x), y' \in T(y).$ It is called \textit{maximal monotone} if in addition its graph $\{(x,x'): x \in \mathcal{H}, x' \in T(x)\}$ is not properly contained in the graph of any monotone operator on $\mathcal{H}.$

\textbf{Remark 2} \cite{subdiffmaximal}: The subdifferential $\partial f$ of a proper closed convex function $f(x)$ is maximal monotone.

\textbf{Remark 3} \cite{rudinfunctionalbook}: In a finite dimensional space, weak convergence implies strong convergence.

\textbf{Lemma 1} [Cor. 2.6]\cite{forwardreflectedbackward}: Let $\Phi: \mathcal{H} \rightrightarrows \mathcal{H}$ be maximally monotone, let $H:\mathcal{H} \longrightarrow \mathcal{H}$ be monotone and $L$-Lipschitz, and suppose that $(\Phi+H)^{-1}(\textbf{0}) \neq \emptyset.$ Choose $\eta \in (0,\frac{1}{2L})$. Given $x_{0},x_{-1} \in \mathcal{H}$, define the sequence $\{x_{k}\}$ according to
\begin{equation}\label{algorithm1}
x_{k+1}=J_{\eta \Phi}(x_{k}-2 \eta H(x_{k})+\eta H(x_{k-1})), \quad{} \forall k \in \mathbb{N}.
\end{equation}
Then $\{x_{k}\}_{-1}^{\infty}$ converges weakly to a point contained in $(\Phi+H)^{-1}(\textbf{0}).$

\textbf{Lemma 2} \cite{threedec}: Let $M \in \Re^{m \times m}$. Then $\Vert M \Vert_{2} \leq \sqrt{\Vert M \Vert_{1} \Vert M \Vert_{\infty}}$.

\textbf{Definition 3:} $B:\mathcal{H} \longrightarrow \mathcal{H}$ is $\beta$-\textit{cocoercive}, $\beta >0,$ if
$$\langle x-y,B(x)-B(y) \rangle \geq \beta \Vert B(x)-B(y) \Vert^{2}, \quad{} \forall x,y \in \mathcal{H}.$$ 

\textbf{Preposition 1} \cite[Preps. 6.4.1 and 6.4.2]{bertsekasconvexbook}: Let the optimal value $f^{*}$ of the following optimization be finite 

\begin{equation}\label{bertsekasbookoptimization}
\begin{aligned}
& \underset{x}{\text{min}}
& & f(x) \\
& \text{subject to}
& & x \in C \\
&&& e_{i}'x-d_{i}=0, \quad{} i=1, \hdots, m, \\
&&& a_{j}'x-b_{j} \leq 0, \quad{} j=1, \hdots, r, 
\end{aligned}
\end{equation}
where $C \subseteq \Re^{n}$ is a nonempty convex set, and the cost function $f:C \longrightarrow \Re$ is convex. If
\begin{align*}
F :=\{x \in C: e_{i}'x-d_{i}=0, \quad{} i=1, \hdots, m, \quad{}
   a_{j}'x-b_{j} \leq 0, \quad{} j=1, \hdots, r \}
\end{align*}
contains a relative interior point of $C$, then the set of geometric multipliers is nonempty.

\textbf{Preposition 2} \cite[Prep. 6.1.2]{bertsekasconvexbook}: Assume that the primal problem (\ref{bertsekasbookoptimization}) has at least one optimal solution $x^{*}$. Then the set of Lagrange multipliers associated with $x^{*}$ and the set of geometric multipliers coincide.

\section{Centralized Optimization}

\subsection{Two-Block Optimization}

Consider optimization (\ref{1}). The \textit{augmented Lagrangian} \cite{boydADMM} for (\ref{1}) is

\begin{align}
\mathcal{L}_{\rho}(x,z,y):= f(x)+g(z)+y^{T} (Ax+Bz-c) +\frac{\rho}{2} \Vert Ax+Bz-c \Vert_{2}^{2} \label{augmentedlagrangian}
\end{align}
where $\rho >0$ is called the \textit{penalty parameter} \cite{boydADMM}, and $y$ is the \textit{Lagrange multiplier} \cite{boydconvexbook} associated with the equality constraint in (\ref{1}). Note that $\mathcal{L}_{0}$ is the standard Lagrangian for the problem (see \cite{boydconvexbook}).

Now we impose the following assumptions on (\ref{1}).

\textbf{Assumption 1:} $f:\Re^{n} \longrightarrow \Re \cup \{+\infty\}$ and $g:\Re^{n} \longrightarrow \Re \cup \{+\infty\}$ are proper, closed, and convex.

Assumption 1 ensures that the single-valued operators $prox_{f}$ and $prox_{g}$ exist (see Remark 1).

\textbf{Assumption 2:} The standard Lagrangian $\mathcal{L}_{0}$ has a saddle point, namely
\begin{align}
 \max_{y \in \Re^{p}} \quad{}  \min_{x \in \Re^{n},z \in \Re^{m}} \mathcal{L}_{0}  (x,z,y)= 
\min_{x \in \Re^{n},z \in \Re^{m}} \quad{}  \max_{y \in \Re^{p}} \quad{} \mathcal{L}_{0}  (x,z,y) \label{saddlepointproblem}
\end{align}
has a solution.

Assumption 2 ensures that the inclusion problem below derived from first order optimality condition of the (augmented) Lagrangian (\ref{augmentedlagrangian}) has a solution. 

Now we give the main theorem in this subsection.

\textbf{Theorem 1:} Consider optimization (\ref{1}) with Assumptions 1 and 2. Let $ \rho \in [0,+\infty)$ and $\eta \in (0,\frac{1}{2L})$ where $\Vert M_{\rho} \Vert_{2} \leq L$ and $M_{\rho}$ is defined in (\ref{matrixM}). Then starting from any initial points $x_{0},x_{-1},z_{0},z_{-1},y_{0},y_{-1}$, the sequences $\{x_{k}\}_{-1}^{\infty}$ and $\{z_{k}\}_{-1}^{\infty}$ generated by Algorithm 1 converges to a solution of (\ref{1}).

\begin{algorithm}
	\caption{Parallel Alternating Direction Primal Dual (PADPD) algorithm}
	\begin{align*}
	\hat{x}_{k} &=x_{k}-2 \eta A^{T} y_{k}-2 \eta \rho A^{T} A x_{k}-2 \eta \rho A^{T} B z_{k} +\eta A^{T} y_{k-1} 
	+ \eta \rho A^{T} A x_{k-1}+ \eta \rho A^{T} B z_{k-1} \nonumber \\
	&\quad{}+\eta \rho A^{T} c \\
	\hat{z}_{k} &= z_{k}-2 \eta B^{T}y_{k}-2 \eta \rho B^{T}A x_{k}-2 \eta \rho B^{T}B z_{k} + \eta B^{T}y_{k-1}+ \eta \rho B^{T}A x_{k-1}+ \eta \rho B^{T}B z_{k-1} \nonumber \\
	&\quad{}+\eta \rho B^{T}c \\
	x_{k+1}&=\displaystyle \arg\min_{u \in \Re^{n}} (\eta f(u)+\frac{1}{2} \Vert u-\hat{x}_{k} \Vert_{2}^{2}) \\
	z_{k+1}&=\displaystyle \arg\min_{r \in \Re^{m}} (\eta g(r)+\frac{1}{2} \Vert r-\hat{z}_{k} \Vert_{2}^{2}) \\
	y_{k+1}&= y_{k}+2 \eta A x_{k}+2 \eta B z_{k}-\eta Ax_{k-1}-\eta Bz_{k-1}-\eta c 
	\end{align*}
\end{algorithm}

\textit{Proof:} Through the first order optimality condition of (\ref{saddlepointproblem}), the saddle point problem (\ref{saddlepointproblem}) can be formulated as the following inclusion problem:

\textit{Find $col\{x^{*},z^{*},y^{*}\}$ such that}
\begin{align}
\textbf{0}_{n} &\in \partial f(x^{*})+A^{T}y^{*}+\rho A^{T} A x^{*}+ \rho A^{T} B z^{*}- \rho A^{T} c \label{inclusion1} \\
\textbf{0}_{m} &\in \partial g(z^{*})+B^{T}y^{*}+\rho B^{T} Ax^{*}+\rho B^{T}Bz^{*}-\rho B^{T}c \label{inclusion2} \\
\textbf{0}_{p} &= -(Ax^{*}+Bz^{*}-c). \label{inclusion3}
\end{align}

Now we consider the Hilbert space $\mathcal{H}=(\Re^{n+m+p}, \Vert . \Vert_{2})$. The inclusion (\ref{inclusion1})-(\ref{inclusion3}) can be rewritten as
\begin{equation}\label{iiiiinclusion}
\textbf{0}_{n+m+p} \in \Phi (\Pi)+H(\Pi)
\end{equation}
where
\begin{align}
\Pi &:=col\{x,z,y\} \\
\Phi (\Pi) &:=col\{ \partial f(x), \partial g(z),\textbf{0}_{p}\} \label{operatorphi} \\
H(\Pi) &:= M_{\rho} \Pi+V_{\rho} \label{operatorH}
\end{align}
in which
\begin{equation}\label{matrixM}
M_{\rho}:=\begin{pmatrix} 
\rho A^{T} A & \rho A^{T}B & A^{T} \\
\rho B^{T}A & \rho B^{T}B & B^{T} \\
-A & -B & \textbf{0}_{p \times p}
\end{pmatrix},
\end{equation}
and
\begin{equation}\label{matrixV}
V_{\rho}:=col\{-\rho A^{T} c, - \rho B^{T} c,c\}.
\end{equation}
From Assumption 1 and Remark 2, we can conclude that the operator $\Phi(\Pi)$ defined in (\ref{operatorphi}) is maximally monotone. Since $\rho \geq 0$, we obtain
\begin{equation}\label{hhhhhh}
(\Pi-\Pi')^{T} (H(\Pi)-H(\Pi')) \geq 0, \quad{} \forall \Pi,\Pi' \in \mathcal{H},
\end{equation}
which implies that the operator $H(\Pi)$ defined in (\ref{operatorH}) is monotone. Since we have not imposed any assumptions on matrices $A$ and $B$, simple calculation shows that the operator $H(\Pi)$ is not cocoercive (see Definition 3). For example, when $\rho=0$, the matrix $M_{0}$ defined in (\ref{matrixM}) is skew symmetric; consequently,
\begin{equation}\label{dddddddddd}
(\Pi-\Pi')^{T} (H(\Pi)-H(\Pi')) = 0, \quad{} \forall \Pi,\Pi' \in \mathcal{H},
\end{equation}
which implies that $H(\Pi)$ is not cocoercive. It is obvious that $H(\Pi)$ defined in (\ref{operatorH}) is $L$-Lipschitz where $\Vert M_{\rho} \Vert_{2} \leq L$.

\textbf{Remark 4:} It seems that applying \textit{forward-backward splitting method} \cite{forwardbackwardfirstalgorithm} to (\ref{iiiiinclusion}) results in a parallel algorithm. However, the forward-backward method \cite{forwardbackwardfirstalgorithm} requires $H(\Pi)$ to be cocoercive (see Definition 3). We have shown above that $H(\Pi)$ is \textit{not} cocoercive. Therefore, we need to apply modified forward-backward method which does not need cocoercivity assumption such as \cite{modifiedforwardbackwardTseng}. Nevertheless, applying \textit{forward-backward-forward method} in \cite{modifiedforwardbackwardTseng} does not yield a parallel algorithm (see \cite[Sec. 4]{modifiedforwardbackwardTseng} for details). \textit{Here, we show that applying \textit{forward-reflected-backward method} without requiring cocoercivity assumption in \cite{forwardreflectedbackward} results in a fully parallel algorithm.}

\textbf{Remark 5:} One may use Lemma 2 to choose $$L=max\{\Vert \begin{pmatrix} 
\rho A^{T} A \\
\rho B^{T}A  \\
-A 
\end{pmatrix} \Vert_{1},\Vert \begin{pmatrix} 
\rho A^{T}B \\
\rho B^{T}B  \\
-B 
\end{pmatrix} \Vert_{1},\Vert [A,B] \Vert_{\infty} \}.$$

Since the saddle-point problem (\ref{saddlepointproblem}) has a solution, the inclusion problem (\ref{inclusion1})-(\ref{inclusion3}) has a solution, i.e., $(\Phi+H)^{-1}(\textbf{0}) \neq \emptyset.$ Therefore, the conditions of Lemma 1 are satisfied, and we can apply Algorithm (\ref{algorithm1}) given $\Pi_{0},\Pi_{-1} \in \mathcal{H},$ i.e.,
\begin{equation}\label{alggggggg}
\Pi_{k+1}=J_{\eta \Phi}(\Pi_{k}-2 \eta H(\Pi_{k})+\eta H(\Pi_{k-1})).
\end{equation}
The sequence $\{\Pi_{k}\}$ generated by (\ref{alggggggg}) converges strongly to a point in $(\Phi+H)^{-1}(\textbf{0})$ since $\mathcal{H}$ is finite dimensional (see Lemma 1 and Remark 3). We obtain
\begin{equation}\label{bbbbbbbbbbb}
\Pi_{k}-2 \eta H(\Pi_{k})+\eta H(\Pi_{k-1})=col\{\Theta_{1,k},\Theta_{2,k},\Theta_{3,k}\}
\end{equation}
where
\begin{align}
\Theta_{1,k} &:=x_{k}-2 \eta A^{T} y_{k}-2 \eta \rho A^{T} A x_{k}-2 \eta \rho A^{T} B z_{k} +\eta A^{T} y_{k-1} 
+ \eta \rho A^{T} A x_{k-1}+ \eta \rho A^{T} B z_{k-1} \nonumber \\
&\quad{}+\eta \rho A^{T} c \\
\Theta_{2,k} &:= z_{k}-2 \eta B^{T}y_{k}-2 \eta \rho B^{T}A x_{k}-2 \eta \rho B^{T}B z_{k} + \eta B^{T}y_{k-1}+ \eta \rho B^{T}A x_{k-1}+ \eta \rho B^{T}B z_{k-1} \nonumber \\
&\quad{}+\eta \rho B^{T}c \\
\Theta_{3,k} &:= y_{k}+2 \eta A x_{k}+2 \eta B z_{k}-\eta c -\eta Ax_{k-1}-\eta Bz_{k-1}.
\end{align} 
We also have
\begin{equation}\label{ppppppppppppp}
J_{\eta \Phi}(\Pi)=\begin{pmatrix} 
\displaystyle \arg\min_{u \in \Re^{n}} (\eta f(u)+\frac{1}{2} \Vert u-x \Vert_{2}^{2}) \\
 \displaystyle \arg\min_{r \in \Re^{m}} (\eta g(r)+\frac{1}{2} \Vert r-z \Vert_{2}^{2}) \\
 \textbf{0}_{p} 
\end{pmatrix}
\end{equation}
which exists by Assumption 1. From (\ref{alggggggg})-(\ref{ppppppppppppp}) we obtain Algorithm 1. We call Algorithm 1 \textit{Parallel Alternating Direction Primal Dual} (PADPD) algorithm since primal variables are updated in alternating fashion, and updates of all primal and dual variables are performed in parallel. This completes the proof of Theorem 1.

A particular interest is for standard Lagrangian $\mathcal{L}_{0}$ defined in (\ref{augmentedlagrangian}). In this case, we have the following corollary.

\textbf{Corollary 1:} Consider optimization (\ref{1}) with Assumptions 1 and 2. Let $\eta \in (0,\frac{1}{2L})$ where $\Vert M_{0} \Vert_{2} \leq L$, and $M_{0}$ is defined in (\ref{matrixM}). Then starting from any initial points $x_{0},x_{-1},z_{0},z_{-1},y_{0},y_{-1}$, the sequences $\{x_{k}\}_{-1}^{\infty}$ and $\{z_{k}\}_{-1}^{\infty}$ generated by Algorithm 2 converges to a solution of (\ref{1}).

\begin{algorithm}
	\caption{}
	\begin{align*}
	\hat{x}_{k} &= x_{k}-2 \eta A^{T} y_{k}+\eta A^{T} y_{k-1}  \nonumber \\
	\hat{z}_{k} &= z_{k}-2 \eta B^{T}y_{k} + \eta B^{T}y_{k-1} \\
	x_{k+1}&=\displaystyle \arg\min_{u \in \Re^{n}} (\eta f(u)+\frac{1}{2} \Vert u-\hat{x}_{k} \Vert_{2}^{2}) \\
	z_{k+1}&=\displaystyle \arg\min_{r \in \Re^{m}} (\eta g(r)+\frac{1}{2} \Vert r-\hat{z}_{k} \Vert_{2}^{2}) \\
	y_{k+1}&= y_{k}+2 \eta A x_{k}+2 \eta B z_{k} -\eta Ax_{k-1}-\eta Bz_{k-1}-\eta c
	\end{align*}
\end{algorithm}

\subsection{Extension to Multi-Block Optimization}

\begin{equation}\label{multiblockoptim}
\begin{aligned}
& \underset{x_{i}, i=1, \hdots,q}{\text{min}}
& & \sum_{i=1}^{q} f_{i}(x_{i}) \\
& \text{subject to}
& & \sum_{i=1}^{q} A_{i} x_{i}=c \\
\end{aligned}
\end{equation}
where $x_{i} \in \Re^{n_{i}}, n_{i} \in \mathbb{N},$ are decision variables, $f_{i}$ are convex functions, and $A_{i} \in \Re^{p \times n_{i}}, i=1, \hdots, q, q \in \mathbb{N}$.

The \textit{augmented Lagrangian} \cite{boydADMM} for (\ref{multiblockoptim}) is
\begin{align}
\mathcal{L}_{\rho}(x_{1},\hdots,x_{q},y):=\sum_{i=1}^{q} f_{i}(x_{i}) +y^{T} (\sum_{i=1}^{q} A_{i} x_{i}-c) +\frac{\rho}{2} \Vert \sum_{i=1}^{q} A_{i} x_{i}-c \Vert_{2}^{2} \label{augmlagrangmultiblock}
\end{align}
where $\rho >0$ is the \textit{penalty parameter} \cite{boydADMM}, and $y$ is the \textit{Lagrange multiplier} \cite{boydconvexbook} associated with the equality constraint in (\ref{multiblockoptim}). $\mathcal{L}_{0}$ is the standard Lagrangian for the problem (see \cite{boydconvexbook}).

We impose the following assumptions on (\ref{multiblockoptim}).

\textbf{Assumption 3:} $f_{i}:\Re^{n_{i}} \longrightarrow \Re \cup \{+\infty\}, i=1, \hdots,q,$ are proper, closed, and convex.

Assumption 3 ensures that the single-valued operators $prox_{f_{i}}, i=1, \hdots,m,$ exist (see Remark 1).

\textbf{Assumption 4:} The standard Lagrangian $\mathcal{L}_{0}$ has a saddle point, i.e.,
\begin{align}
 \max_{y \in \Re^{p}} \quad{} \quad{}  \min_{x_{i} \in \Re^{n_{i}},i=1,\hdots,q} \quad{} \mathcal{L}_{0}  (x_{1},\hdots,x_{q},y)= \min_{x_{i} \in \Re^{n_{i}},i=1,\hdots,q} \quad{} \quad{}  \max_{y \in \Re^{p}} \quad{} \quad{}  \mathcal{L}_{0}  (x_{1},\hdots,x_{q},y) \label{saddlepointproblemmulti}
\end{align}
has a solution.

Assumption 4 ensures that the inclusion problem below derived from first order optimality condition of the (augmented) Lagrangian (\ref{augmlagrangmultiblock}) has a solution.

\begin{algorithm}
	\caption{Parallel Alternating Direction Primal Dual (PADPD) Algorithm for Multi-Block Optimization (\ref{multiblockoptim})}
	\begin{align*}
	\hat{x}_{1,k} &=x_{1,k}-2 \eta A_{1}^{T} y_{k}-2 \eta \rho \sum_{j=1}^{q} A_{1}^{T} A_{j} x_{j,k} +\eta A_{1}^{T} y_{k-1}+\eta \rho \sum_{j=1}^{q} A_{1}^{T} A_{j} x_{j,k-1}+\rho A_{1}^{T} c \\
	&\quad{} \vdots \nonumber \\
	\hat{x}_{i,k} &=x_{i,k}-2 \eta A_{i}^{T} y_{k}-2 \eta \rho \sum_{j=1}^{q} A_{i}^{T} A_{j} x_{j,k} +\eta A_{i}^{T} y_{k-1}+\eta \rho \sum_{j=1}^{q} A_{i}^{T} A_{j} x_{j,k-1}+\rho A_{i}^{T} c \\
	&\quad{} \vdots \nonumber \\
	\hat{x}_{q,k} &=x_{q,k}-2 \eta A_{q}^{T} y_{k}-2 \eta \rho \sum_{j=1}^{q} A_{q}^{T} A_{j} x_{j,k} +\eta A_{q}^{T} y_{k-1}+\eta \rho \sum_{j=1}^{q} A_{q}^{T} A_{j} x_{j,k-1}+\rho A_{q}^{T} c \\
	x_{1,k+1}&=\displaystyle \arg\min_{u_{1} \in \Re^{n_{1}}} (\eta f_{1}(u_{1})+\frac{1}{2} \Vert u_{1}-\hat{x}_{1,k} \Vert_{2}^{2})  \\
	&\quad{} \vdots \nonumber \\
	x_{i,k+1}&=\displaystyle \arg\min_{u_{i} \in \Re^{n_{i}}} (\eta f_{i}(u_{i})+\frac{1}{2} \Vert u_{i}-\hat{x}_{i,k} \Vert_{2}^{2})  \\
	&\quad{} \vdots \nonumber \\
	x_{q,k+1}&=\displaystyle \arg\min_{u_{q} \in \Re^{n_{q}}} (\eta f_{q}(u_{q})+\frac{1}{2} \Vert u_{q}-\hat{x}_{q,k} \Vert_{2}^{2})  \\
	y_{k+1}&=y_{k}+2 \eta \sum_{i=1}^{q} A_{i}x_{i,k}-\eta \sum_{i=1}^{q} A_{i}x_{i,k-1} -\eta c
	\end{align*}
\end{algorithm}

\textbf{Theorem 2:} Consider optimization (\ref{multiblockoptim}) with Assumptions 3 and 4. Let $ \rho \in [0,+\infty)$ and $\eta \in (0,\frac{1}{2 \tilde{L}})$ where $\Vert \tilde{M}_{\rho} \Vert_{2} \leq \tilde{L}$ and $\tilde{M}_{\rho}$ is defined in (\ref{matrixMmulti}). Then starting from any initial points $x_{1,0},x_{1,-1},\hdots,x_{q,0},x_{q,-1},y_{0},y_{-1}$, the sequences $\{x_{i,k}\}_{-1}^{\infty}, i=1,\hdots,q,$ generated by Algorithm 3 converge to a solution of (\ref{multiblockoptim}).

\textit{Proof:} Through the first optimality condition of (\ref{saddlepointproblemmulti}), the saddle point problem (\ref{saddlepointproblemmulti}) can be formulated as the following inclusion problem:

\textit{Find $col\{x_{1}^{*},\hdots,x_{q}^{*},y^{*}\}$ such that} 
\begin{align}
\textbf{0}_{n_{1}} &\in \partial f_{1}(x_{1}^{*}) +A_{1}^{T} y^{*} +\rho \sum_{j=1}^{q} A_{1}^{T} A_{j} x_{j}^{*}-\rho A_{1}^{T} c \label{inclusionmulti1} \\
\textbf{0}_{n_{2}} &\in \partial f_{2}(x_{2}^{*}) +A_{2}^{T} y^{*} +\rho \sum_{j=1}^{q} A_{2}^{T} A_{j} x_{j}^{*}-\rho A_{2}^{T} c \label{inclusionmulti2} \\
&\vdots \nonumber \\
\textbf{0}_{n_{i}} &\in \partial f_{i}(x_{i}^{*}) +A_{i}^{T} y^{*} +\rho \sum_{j=1}^{q} A_{i}^{T} A_{j} x_{j}^{*}-\rho A_{i}^{T} c \label{inclusionmulti3} \\
&\vdots \nonumber \\
\textbf{0}_{p} &= -(\sum_{i=1}^{q} A_{i} x_{i}-c). \label{inclusionmulti4}
\end{align}

Now we consider the Hilbert space $\tilde{\mathcal{H}}=(\Re^{p+\sum_{i=1}^{q} n_{i}}, \Vert . \Vert_{2})$. The inclusion (\ref{inclusionmulti1})-(\ref{inclusionmulti4}) can be rewritten as

\begin{equation}
\textbf{0}_{\Re^{p+\sum_{i=1}^{q} n_{i}}} \in \tilde{\Phi} (\tilde{\Pi})+\tilde{H}(\tilde{\Pi})
\end{equation}
where
\begin{align}
\tilde{\Pi} &:=col\{x_{1},\hdots,x_{m},y\} \\
\tilde{\Phi} (\tilde{\Pi}) &:=col\{ \partial f_{1}(x_{1}), \hdots, \partial f_{q}(x_{q}),\textbf{0}_{p}\} \label{operatorphimulti} \\
\tilde{H}(\tilde{\Pi}) &:= \tilde{M}_{\rho} \tilde{\Pi}+\tilde{V}_{\rho} \label{operatorHmulti}
\end{align}
in which
\begin{equation}\label{matrixMmulti}
\tilde{M}_{\rho}:=\begin{pmatrix} 
\rho A^{T}_{1} A_{1} & \rho A^{T}_{1}A_{2} & \hdots & \rho A^{T}_{1}A_{q} & A^{T}_{1} \\
\rho A^{T}_{2} A_{1} & \rho A^{T}_{2}A_{2} & \hdots & \rho A^{T}_{2}A_{q} & A^{T}_{2} \\
\vdots & \vdots & \ddots & \vdots & \vdots \\
\rho A^{T}_{q} A_{1} & \rho A^{T}_{q}A_{2} & \hdots & \rho A^{T}_{q}A_{q} & A^{T}_{q} \\
-A_{1} & -A_{2} & \hdots & -A_{q} & \textbf{0}_{p \times p}
\end{pmatrix},
\end{equation}
and
\begin{equation}\label{matrixVmulti}
\tilde{V}_{\rho}:=col\{-\rho A^{T}_{1} c, - \rho A^{T}_{2} c, \hdots,- \rho A^{T}_{q} c, c\}.
\end{equation}

From Assumption 3 and Remark 2, we conclude that the operator $\tilde{\Phi} (\tilde{\Pi})$ defined in (\ref{operatorphimulti}) is maximally monotone. Since $\rho \geq 0,$ we obtain
\begin{equation}\label{hhhhhhmulti}
(\tilde{\Pi}-\tilde{\Pi}')^{T} (\tilde{H}(\tilde{\Pi})-\tilde{H}(\tilde{\Pi}')) \geq 0, \quad{} \forall \tilde{\Pi},\tilde{\Pi}' \in \tilde{\mathcal{H}},
\end{equation}
which implies that the operator $\tilde{H}(\tilde{\Pi})$ defined in (\ref{operatorHmulti}) is monotone. Since we have not imposed any assumptions on matrices $A_{i},i=1\hdots,m,$ one can show that the operator $\tilde{H}(\tilde{\Pi})$ is not cocoercive (see Definition 3). It is obvious that $\tilde{H}(\tilde{\Pi})$ defined in (\ref{operatorHmulti}) is $\tilde{L}$-Lipschitz where $\Vert \tilde{M}_{\rho} \Vert_{2} \leq \tilde{L}$.

\textbf{Remark 6:} One may use Lemma 2 to choose 
\begin{align*}
\tilde{L}&=max\{\Vert \begin{pmatrix} 
\rho A^{T}_{1} A_{1} \\
\rho A^{T}_{2} A_{1}  \\
\vdots \\
\rho A^{T}_{q} A_{1}  \\
-A_{1}
\end{pmatrix} \Vert_{1},\Vert \begin{pmatrix} 
\rho A^{T}_{1} A_{2} \\
\rho A^{T}_{2} A_{2}  \\
\vdots \\
\rho A^{T}_{q} A_{2}  \\
-A_{2}
\end{pmatrix} \Vert_{1}, \hdots, \Vert \begin{pmatrix} 
\rho A^{T}_{1} A_{q} \\
\rho A^{T}_{2} A_{q}  \\
\vdots \\
\rho A^{T}_{q} A_{q}  \\
-A_{q}
\end{pmatrix} \Vert_{1}, \Vert [A_{1},A_{2},\hdots,A_{q}] \Vert_{\infty} \}.
\end{align*}

Since the saddle-point problem (\ref{saddlepointproblemmulti}) has a solution, the inclusion problem (\ref{inclusionmulti1})-(\ref{inclusionmulti4}) has a solution, i.e., $(\tilde{\Phi}+\tilde{H})^{-1}(\textbf{0}) \neq \emptyset.$ Therefore, the conditions of Lemma 1 are satisfied, and we can apply Algorithm (\ref{algorithm1}) given $\tilde{\Pi}_{0},\tilde{\Pi}_{-1} \in \tilde{\mathcal{H}},$ i.e.,
\begin{equation}\label{algggggggmulti}
\tilde{\Pi}_{k+1}=J_{\eta \tilde{\Phi}}(\tilde{\Pi}_{k}-2 \eta \tilde{H}(\tilde{\Pi}_{k})+\eta \tilde{H}( \tilde{\Pi}_{k-1})).
\end{equation}
The sequence $\{\tilde{\Pi}_{k}\}$ generated by (\ref{algggggggmulti}) converges strongly to a point in $(\tilde{\Phi}+\tilde{H})^{-1}(\textbf{0})$ since $\tilde{\mathcal{H}}$ is finite dimensional (see Lemma 1 and Remark 3). We obtain
\begin{align}
&\tilde{\Pi}_{k}-2 \eta \tilde{H}(\tilde{\Pi}_{k})+\eta \tilde{H}(\tilde{\Pi}_{k-1})= \nonumber \\
&col\{\tilde{\Theta}_{1,k},\hdots,\tilde{\Theta}_{i,k},\hdots, \tilde{\Theta}_{q+1,k}\} \label{bbbbbbbbbbbmulti}
\end{align}
where
\begin{align}
\tilde{\Theta}_{1,k} &:=x_{1,k}-2 \eta A_{1}^{T} y_{k}-2 \eta \rho \sum_{j=1}^{q} A_{1}^{T} A_{j} x_{j,k} +\eta A_{1}^{T} y_{k-1}+\eta \rho \sum_{j=1}^{q} A_{1}^{T} A_{j} x_{j,k-1}+\rho A_{1}^{T} c \\
&\quad{} \vdots \nonumber \\
\tilde{\Theta}_{i,k} &:=x_{i,k}-2 \eta A_{i}^{T} y_{k}-2 \eta \rho \sum_{j=1}^{q} A_{i}^{T} A_{j} x_{j,k} +\eta A_{i}^{T} y_{k-1}+\eta \rho \sum_{j=1}^{q} A_{i}^{T} A_{j} x_{j,k-1}+\rho A_{i}^{T} c \\
&\quad{} \vdots \nonumber \\
\tilde{\Theta}_{q+1,k} &:=y_{k}+2 \eta \sum_{i=1}^{q} A_{i}x_{i,k}-\eta \sum_{i=1}^{q} A_{i}x_{i,k-1} -\eta c.
\end{align}
We have that
\begin{equation}\label{pppppppppppppmulti}
J_{\eta \tilde{\Phi}}(\tilde{\Pi})=\begin{pmatrix} 
\displaystyle \arg\min_{u_{1} \in \Re^{n_{1}}} (\eta f_{1}(u_{1})+\frac{1}{2} \Vert u_{1}-x_{1} \Vert_{2}^{2}) \\
\vdots \\
\displaystyle \arg\min_{u_{q} \in \Re^{n_{q}}} (\eta f_{q}(u_{q})+\frac{1}{2} \Vert u_{q}-x_{q} \Vert_{2}^{2}) \\
\textbf{0}_{p} 
\end{pmatrix}
\end{equation}
which exists by Assumption 3. From (\ref{algggggggmulti})-(\ref{pppppppppppppmulti}) we obtain Algorithm 3. Thus the proof of Theorem 2 is complete.

\begin{algorithm}
	\caption{}
	\begin{align*}
	\hat{x}_{1,k} &= x_{1,k}-2 \eta A^{T}_{1} y_{k}+\eta A^{T}_{1} y_{k-1}  \nonumber \\
	&\quad{} \vdots \\
	\hat{x}_{i,k} &= x_{i,k}-2 \eta A^{T}_{i} y_{k}+\eta A^{T}_{i} y_{k-1}  \nonumber \\
	&\quad{} \vdots \\
	\hat{x}_{q,k} &= x_{q,k}-2 \eta A^{T}_{q} y_{k}+\eta A^{T}_{q} y_{k-1}  \nonumber \\
	x_{1,k+1}&=\displaystyle \arg\min_{u_{1} \in \Re^{n_{1}}} (\eta f_{1}(u_{1})+\frac{1}{2} \Vert u_{1}-\hat{x}_{1,k} \Vert_{2}^{2}) \\
	&\quad{} \vdots \\
	x_{i,k+1}&=\displaystyle \arg\min_{u_{i} \in \Re^{n_{i}}} (\eta f_{i}(u_{i})+\frac{1}{2} \Vert u_{i}-\hat{x}_{i,k} \Vert_{2}^{2}) \\
	&\quad{} \vdots \\
	x_{q,k+1}&=\displaystyle \arg\min_{u_{q} \in \Re^{n_{q}}} (\eta f_{q}(u_{q})+\frac{1}{2} \Vert u_{q}-\hat{x}_{q,k} \Vert_{2}^{2}) \\
	y_{k+1}&=y_{k}+2 \eta \sum_{i=1}^{q} A_{i}x_{i,k}-\eta \sum_{i=1}^{q} A_{i}x_{i,k-1} -\eta c
	\end{align*}
\end{algorithm}

A particular interest is for standard Lagrangian $\mathcal{L}_{0}$ defined in (\ref{augmlagrangmultiblock}). In this case, we have the following corollary.

\textbf{Corollary 2:} Consider optimization (\ref{multiblockoptim}) with Assumptions 3 and 4. Let $\eta \in (0,\frac{1}{2 \tilde{L}})$ where $\Vert \tilde{M}_{0} \Vert_{2} \leq \tilde{L}$ and $\tilde{M}_{0}$ is defined in (\ref{matrixMmulti}). Then starting from any initial points $x_{1,0},x_{1,-1},\hdots,$ $x_{q,0},x_{q,-1},y_{0},y_{-1}$, the sequences $\{x_{i,k}\}_{-1}^{\infty}, i=1,\hdots,m,$ generated by Algorithm 4 converge to a solution of (\ref{multiblockoptim}).

\section{Distributed Optimization}

Now we consider a network model for the distributed optimization considered below. A network of $m \in \mathbb{N}$ nodes labeled by the set $\mathcal{V}=\lbrace 1,2,...,m \rbrace $ is considered. The topology of the interconnections among nodes is fixed $\mathcal{G}=(\mathcal{V},\mathcal{E})$ where $\mathcal{E}$ is the \textit{unordered} edge set $\mathcal{E} \subseteq \mathcal{V} \times \mathcal{V}$. We write $\mathcal{N}_{i}$ for the labels of agent $i$'s neighbors. We define the weighted graph matrix $\mathcal{W}=[\mathcal{W}_{ij}]$  with $\mathcal{W}_{ij}=a_{ij}$ for $j \in \mathcal{N}_{i} \cup \{ i \}$, and $\mathcal{W}_{ij}=0$ otherwise, where $a_{ij}>0$ is the scalar constant weight that agent $i$ assigns to the information $x_{j}$ received from agent $j$.

Now we define the distributed optimization problem as follows: for each node $i \in \mathcal{V}$, we associate a private convex cost function $f_{i}:\Re^{n} \longrightarrow \Re$ which is known to node $i$. The objective of each agent is to collaboratively seek the solution of the following optimization problem using local information exchange with the neighbors:
$$\underset{s}\min \sum_{i=1}^{m} f_{i}(s)$$
where $s \in \Re^{n}$. We assume that there is no communication delay or noise in delivering a message from agent $j$ to agent $i$.

The full formulation of the above distributed optimization problem is as follows:

\begin{equation}\label{11}
\begin{aligned}
& \underset{x}{\text{min}}
& & f(x):=\sum_{i=1}^{m} f_{i}(x_{i}) \\
& \text{subject to}
& & x_{1}=x_{2}=...=x_{m}
\end{aligned}
\end{equation}
where $x=[x_{1}^{T},...,x_{m}^{T}]^{T}, x_{i} \in \Re^{n}, i=1,2,...,m$, $f_{i}:\Re^{n} \longrightarrow \Re$ is a private cost function known to node $i$, and the constraint is achieved through interactions with neighbors.

Now we impose the following assumptions on (\ref{11}).

\textbf{Assumption 5:} The solution set $\mathcal{X}^{*}$ of optimization (\ref{11}) is nonempty.

\textbf{Assumption 6:} The cost functions $f_{i}:\Re^{n} \longrightarrow \Re, i=1, \hdots, m,$ are convex.

Assumptions 5-6 imply that the solution set of (\ref{11}) is nonempty, closed, and convex. Under Assumption 6, strong duality \cite{boydconvexbook} holds for (\ref{11}) and its Lagrange's dual problem (by Weak Slater's condition). Moreover, under Assumptions 5-6, the conditions of Preposition 1 are satisfied, and we can guarantee the existence of a Lagrange multiplier associated with $s^{*} \in \mathcal{X}^{*}$ by Preposition 2 since (\ref{11}) is a special case of (\ref{bertsekasbookoptimization}). Note that Assumption 1 is satisfied, for example, if the cost functions $f_{i}$ satisfy some growth conditions.

Now we impose the following assumptions on the weighted matrix of the graph.

\textbf{Assumption 7:} $\mathcal{W}=\mathcal{W}^{T}$ and $\mathcal{W} \textbf{1}_{m}=\textbf{1}_{m}.$

Assumption 7 implies that the links in the graph are undirected, and the weighted matrix of the graph is doubly stochastic.

\textbf{Assumption 8} \cite{alavianiACC2017}-\cite{alavianiTAC}: The graph is connected, i.e., $\lambda_{2} (I_{m}-\mathcal{W}) >0.$ 

Assumption 8 ensures that the information sent from each node will be finally obtained by every other node through a path.

Now, a solution to (\ref{11}) under Assumptions 7-8 can be obtained by solving the following problem:

\begin{equation}\label{12}
\begin{aligned}
& \underset{x}{\text{min}}
& & f(x):=\sum_{i=1}^{m} f_{i}(x_{i}) \\
& \text{subject to}
& & Wx=x
\end{aligned}
\end{equation}
where $W=\mathcal{W} \otimes I_{n}$. The proof that a solution to (\ref{12}) provides a solution to (\ref{11}) is given in \cite[App. B]{alavianiTAC}. Now we have the following corollary for distributed optimization (\ref{12}).

\textbf{Corollary 3}: Consider distributed optimization (\ref{12}) with Assumptions 5-8. Let $\eta \in (0,\frac{1}{4})$. Then starting from any initial points $x_{i,0},x_{i,-1}, i=1,\hdots,m,$ the sequences $\{x_{i,k}\}_{-1}^{\infty}$ generated by Algorithm 5 converge to a solution of (\ref{11}).

\textbf{Remark 7:} It is obvious that choosing the parameter $\eta$ in Corollary 3 does \textit{not} require structure of the graph or any global information.
\begin{algorithm}
	\caption{}
	\begin{align*}
	\hat{x}_{i,k} &= x_{i,k} +(\eta y_{i,k-1}-2 \eta y_{i,k})-\sum_{j \in \mathcal{N}_{i} \cup \{i\}} \mathcal{W}_{ij} (\eta y_{j,k-1}-2 \eta y_{j,k})  \nonumber \\
	x_{i,k+1}&=\displaystyle \arg\min_{u_{i} \in \Re^{n}} (\eta f_{i}(u_{i})+\frac{1}{2} \Vert u_{i}-\hat{x}_{i,k} \Vert_{2}^{2}) \\
	y_{i,k+1}&= y_{i,k}+ (2 \eta x_{i,k} -\eta x_{i,k-1}) -\sum_{j \in \mathcal{N}_{i} \cup \{i\}} \mathcal{W}_{ij} (2 \eta x_{j,k} -\eta x_{j,k-1})  
	\end{align*}
\end{algorithm}

\textit{Proof of Corollary 3:} It is clear that optimization (\ref{12}) is of the form (\ref{1}). Hence, we can apply Corollary 1 once its conditions are satisfied. By Assumptions 5-6, a Lagrange multiplier associated with a solution of (\ref{12}) exists from Prepositions 1-2. Therefore, Assumptions 1-2 are satisfied. Hence, the conditions of Corollary 1 are satisfied, and we can apply Algorithm 2 for the above distributed optimization problem. By the fact that 
$$\Vert I_{mn}-W \Vert_{\infty} \leq 2 \quad{} \text{and} \quad{}\Vert I_{mn}-W \Vert_{1} \leq 2$$
(see Assumption 7), we obtain from Lemma 2 that
$$\Vert M_{0} \Vert_{2}^{2} \leq \Vert M_{0} \Vert_{1} \Vert M_{0} \Vert_{\infty}= 4 \quad{} \Rightarrow \quad{} L=2.$$
 Algorithm 2 reduces to the following distributed algorithm:
\begin{align*}
\hat{x}_{k} &= x_{k}-2 \eta (I_{mn}-W) y_{k}+\eta (I_{mn}-W) y_{k-1}  \nonumber \\
x_{k+1}&=\displaystyle \arg\min_{u \in \Re^{mn}} (\eta f(u)+\frac{1}{2} \Vert u-\hat{x}_{k} \Vert_{2}^{2}) \\
y_{k+1}&= y_{k}+2 \eta (I_{mn}-W) x_{k} -\eta (I_{mn}-W)x_{k-1}
\end{align*}
The above algorithm is in a compact form and can be  viewed as Algorithm 5 based
on local information for each agent $i$. Thus the proof of Corollary 3 is complete.

\section{Numerical Example}

We give the following example where the \textit{direct} extension algorithm of ADMM is divergent, but Algorithm 3 can converge to a solution of the problem.

\textit{Example 1}: Consider the following three-block optimization \cite[Rem. 3.2]{ADMMdiverges}:

\begin{equation}\label{ADMMdivergentoptim}
\begin{aligned}
& \underset{}{\text{min}}
& & \frac{1}{2} x_{1}^{2} \\
& \text{subject to}
& & \begin{pmatrix} 
1 & 1\\
1 & 1 \\
1 & 1
\end{pmatrix}\begin{pmatrix} 
x_{1} \\
x_{2}
\end{pmatrix}+\begin{pmatrix} 
1 \\
1  \\
2
\end{pmatrix} x_{3}+\begin{pmatrix} 
1 \\
2 \\
2
\end{pmatrix} x_{4}=0. \\
\end{aligned}
\end{equation}
Optimization (\ref{ADMMdivergentoptim}) has a unique optimal solution $x_{1}=x_{2}=x_{3}=x_{4}=0$ for which the optimal value is finite. Consequently, the conditions of Prepositions 1 and 2 are satisfied, and we can guarantee existence of Lagrange multiplier associated with the unique solution. Therefore, Assumption 4 is satisfied. Hence, conditions of Theorem 2 are satisfied.

It has been shown analytically in \cite{ADMMdiverges} that direct extension of ADMM for (\ref{ADMMdivergentoptim}) diverges for \textit{any} $\rho >0$ in augmented Lagrangian (\ref{augmlagrangmultiblock}) and \textit{any} initial conditions. Now we simulate Algorithm 3 where $\rho=1$ and Algorithm 4 (which is Algorithm 3 where $\rho=0$) and show that the sequences generated by the algorithms converge to the solution of (\ref{ADMMdivergentoptim}) by selecting random initial conditions.

Here, we have that
$$A_{1}=\begin{pmatrix} 
1 & 1 \\
1 & 1 \\
1 & 1
\end{pmatrix}, A_{2}=\begin{pmatrix} 
1 \\
1 \\
2
\end{pmatrix}, A_{3}=\begin{pmatrix} 
1 \\
2 \\
2
\end{pmatrix}, c= \textbf{0}_{3},$$
$$ f_{1}(x_{1},x_{2})=\frac{1}{2} x_{1}^{2}, f_{2}(x_{3})=0, f_{3}(x_{4})=0.$$

Hence, Algorithm 3 reduces to the following algorithm:
\begin{align}
\begin{pmatrix} 
x_{1,k+1} \\
x_{2,k+1}
\end{pmatrix} &=\begin{pmatrix} 
\frac{1}{1+\eta} & 0 \\
0 & 1
\end{pmatrix}[\begin{pmatrix} 
1-6 \eta \rho & -6 \eta \rho \\
-6 \eta \rho & 1-6 \eta \rho
\end{pmatrix} \begin{pmatrix} 
x_{1,k} \\
x_{2,k}
\end{pmatrix}  -2 \eta \begin{pmatrix} 
1 & 1 & 1 \\
1 & 1 & 1
\end{pmatrix} y_{k} -2 \eta \rho \begin{pmatrix} 
4 \\
4
\end{pmatrix} x_{3,k} \nonumber \\
&\quad{}-2 \eta \rho \begin{pmatrix} 
5 \\
5
\end{pmatrix} x_{4,k}+\eta \begin{pmatrix} 
1 & 1 & 1 \\
1 & 1 & 1
\end{pmatrix} y_{k-1} +\eta \rho  \begin{pmatrix} 
3 & 3 \\
3 & 3
\end{pmatrix} \begin{pmatrix} 
x_{1,k-1} \\
x_{2,k-1}
\end{pmatrix} +\eta \rho \begin{pmatrix} 
4 \\
4
\end{pmatrix} x_{3,k-1} \nonumber \\
&\quad{}+\eta \rho \begin{pmatrix} 
5\\
5
\end{pmatrix} x_{4,k-1}  ] \label{examplealg1} 
\end{align}

\begin{align}
x_{3,k+1} &=(1-12 \eta \rho) x_{3,k}-2 \eta \begin{pmatrix} 
1 & 1 & 2
\end{pmatrix} y_{k} -2 \eta \rho \begin{pmatrix} 
4 & 4
\end{pmatrix} \begin{pmatrix} 
x_{1,k} \\
x_{2,k}
\end{pmatrix} -14 \eta \rho x_{4,k} \nonumber \\
&\quad{}+\eta \begin{pmatrix} 
1 & 1 & 2
\end{pmatrix} y_{k-1}+\eta \rho \begin{pmatrix} 
4 & 4
\end{pmatrix} \begin{pmatrix} 
x_{1,k-1} \\
x_{2,k-1}
\end{pmatrix} + 6 \eta \rho x_{3,k-1}+7 \eta \rho x_{4,k-1}  \\
x_{4,k+1} &=(1-18 \eta \rho)x_{4,k} 
-2 \eta \begin{pmatrix} 
1 & 2 & 2 
\end{pmatrix} y_{k} -2 \eta \rho \begin{pmatrix} 
5 & 5
\end{pmatrix} \begin{pmatrix} 
x_{1,k} \\
x_{2,k}
\end{pmatrix} -14 \eta \rho x_{3,k} \nonumber \\
&\quad{} +\eta \begin{pmatrix} 
1 & 2 & 2
\end{pmatrix} y_{k-1} +\eta \rho \begin{pmatrix} 
5 & 5
\end{pmatrix} \begin{pmatrix} 
x_{1,k-1} \\
x_{2,k-1}
\end{pmatrix}  + 7 \eta \rho x_{3,k-1}+9 \eta \rho x_{4,k-1} \\
y_{k+1} &= y_{k}+2 \eta \begin{pmatrix} 
1 & 1 \\
1 & 1 \\
1 & 1
\end{pmatrix} \begin{pmatrix} 
x_{1,k} \\
x_{2,k}
\end{pmatrix} +2 \eta \begin{pmatrix} 
1 \\
1 \\
2
\end{pmatrix} x_{3,k} +2 \eta \begin{pmatrix} 
1 \\
2 \\
2
\end{pmatrix} x_{4,k}- \eta \begin{pmatrix} 
1 & 1 \\
1 & 1 \\
1 & 1
\end{pmatrix}\begin{pmatrix} 
x_{1,k-1} \\
x_{2,k-1}
\end{pmatrix} \nonumber \\
&\quad{}- \eta \begin{pmatrix} 
1 \\
1 \\
2
\end{pmatrix} x_{3,k-1}-\eta \begin{pmatrix} 
1 \\
2 \\
2
\end{pmatrix} x_{4,k-1} \label{examplealg2}.
\end{align}

\begin{figure}[thpb]
	\centering
	\includegraphics[scale=0.6]{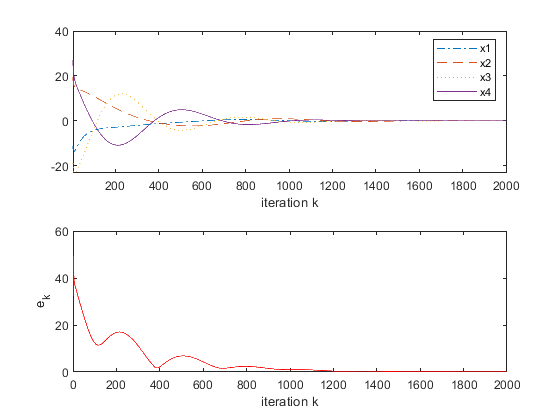}
	\caption{Top: Variables $x_{1,k},x_{2,k},x_{3,k},$ and $x_{4,k}$ of Algorithm (\ref{examplealg1})-(\ref{examplealg2}) where $\rho=1$. The result shows that the variables are approaching the origin. Bottom: The error $e_{k}=\Vert [x_{1,k},x_{2,k},x_{3,k},x_{4,k}] \Vert_{2}$. The result shows that the error is converging to zero.}
	\label{figure1}
\end{figure}

\begin{figure}[thpb]
	\centering
	\includegraphics[scale=0.6]{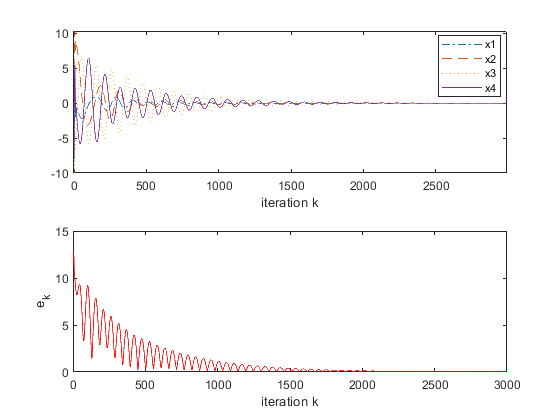}
	\caption{Top: Variables $x_{1,k},x_{2,k},x_{3,k},$ and $x_{4,k}$ of Algorithm (\ref{examplealg1})-(\ref{examplealg2}) where $\rho=0$. The result shows that the variables are approaching the origin. Bottom: The error $e_{k}=\Vert [x_{1,k},x_{2,k},x_{3,k},x_{4,k}] \Vert_{2}$. The result shows that the error is converging to zero.}
	\label{figure2}
\end{figure}

We use MATLAB software for simulation. First, we set $\rho=1$ in Algorithm (\ref{examplealg1})-(\ref{examplealg2}) and calculate $\Vert \tilde{M}_{1} \Vert_{2}=21.3217$. Hence, we select $\eta=\frac{1}{50}$. The result obtained by the algorithm is given in Figure \ref{figure1}. Finally, we set $\rho=0$ in 
Algorithm (\ref{examplealg1})-(\ref{examplealg2})  (to obtain Algorithm 4) and calculate $\Vert \tilde{M}_{0} \Vert_{2}=4.5129$. Hence, we select $\eta=0.1$. The result obtained by the algorithm is given in Figure \ref{figure2}. The results in both figures show that the sequences $\{x_{1,k}\},\{x_{2,k}\},\{x_{3,k}\},$ and $\{x_{4,k}\}$ generated by Algorithms 3 and 4 converge to the solution of (\ref{ADMMdivergentoptim}).

\section{Conclusions and Future Work}

In this paper, we consider a centralized two-block separable optimization for which we derive a fully parallel primal-dual discrete-time algorithm based on monotone operator splitting method. In this algorithm, the primal variables are updated in an alternating fashion like ADMM. However, unlike existing discrete-time algorithms such as MM, ADMM, BiADMM, and PDFP algorithms, that all suffer from sequential updates, all primal and dual variables are updated in parallel. One of advantages of the proposed algorithm is that it can be directly extended to any finite multi-block optimization while preserving its convergence. Then the method is applied to distributed optimization to derive a fully parallel primal dual distributed algorithm. Finally, a numerical example of a three-block optimization is given for which the direct extension of proposed algorithm is shown to converge to a solution, whereas the direct extension of ADMM diverges for any choice of $\rho >0$ and any initial conditions. Rate of convergence of the proposed algorithms remains a topic of future research.

\end{document}